\documentclass[final, preprint]{elsarticle}

\usepackage{fourier}
\usepackage[T1]{fontenc}
\usepackage{graphicx}
\usepackage{xcolor}
\usepackage{float}
\usepackage{tikz}
\usepackage{subcaption}
\usepackage{amsmath,amssymb}
\usepackage{palatino}

\usepackage{url}
\newtheorem{fact}{Fact}
\newtheorem{theorem}{Theorem}
\newtheorem{lemma}[theorem]{Lemma}
\newtheorem{corollary}[theorem]{Corollary}
\newdefinition{rmk}{Remark}
\newproof{proof}{Proof}

\begin{document}

\journal{Information Processing Letters}

%

%

\begin{frontmatter}






\title{New Results on Pairwise Compatibility Graphs}


\author[First]{Sheikh Azizul Hakim\corref{cor1}}
\ead{1705002@ugrad.cse.buet.ac.bd}

\author[Second]{Bishal Basak Papan}
\ead{1505043.bbp@ugrad.cse.buet.ac.bd}

\author[Third]{Md. Saidur Rahman}
\ead{saidurrahman@cse.buet.ac.bd}

\cortext[cor1]{Corresponding author}

\address{Graph Drawing and Information Visualization Laboratory, \\
Department of Computer Science and Engineering, \\
Bangladesh University of Engineering and Technology (BUET), \\
Dhaka-1205, Bangladesh.}

\begin{abstract}
A graph $G=(V,E)$ is called a pairwise compatibility graph (PCG) if there exists an edge-weighted tree $T$ and two non-negative real numbers $d_{min}$ and $d_{max}$ such that each leaf $u$ of $T$ corresponds to a vertex $u \in V$ and there is an edge $(u, v) \in E$ if and only if $d_{min} \leq d_{T}(u, v) \leq d_{max}$, where $d_T(u, v)$ is the sum of the weights of the edges on the unique path from $u$ to $v$ in $T$. The tree $T$ is called the pairwise compatibility tree (PCT) of $G$. It has been proven that not all graphs are PCGs. Thus, it is interesting to know which classes of graphs are PCGs. In this paper, we prove that grid graphs are PCGs. Although there are a necessary condition and a sufficient condition known for a graph being a PCG, there are some classes of graphs that are intermediate to the classes defined by the necessary condition and the sufficient condition. In this paper, we show two examples of graphs that are included in these intermediate classes and prove that they are not PCGs. 
\end{abstract}


\end{frontmatter}

\section{Introduction}

\begin{figure}
\begin{minipage}{0.55\linewidth}
\begin{center}
\begin{tikzpicture}[roundnode/.style={circle, draw=black, inner sep=0pt, minimum size=1.5mm}]

        \node[roundnode, fill=black, label=90:{\footnotesize{$u_3$}}] at (-2.5,0.25) (a3) {};
        \node[roundnode, fill=black, label=180:{\footnotesize{$u_2$}}] at (-3.25,-0.75) (a2) {};
        \node[roundnode, fill=black, label=270:{\footnotesize{$u_1$}}] at (-2.5,-1.75) (a1) {};
        \node[roundnode, fill=white, label=0:{}] at (-2,-0.75) (aa0) {};
        \node[roundnode, fill=white, label=0:{}] at (-1.25,-0.75) (aa1) {};
        \node[roundnode, fill=black, label=90:{\footnotesize{$u_4$}}] at (-1.25,0.25) (a4) {};
        \node[roundnode, fill=white, label=0:{}] at (-0.5,-0.75) (aa2) {};
        \node[roundnode, fill=black, label=90:{\footnotesize{$u_5$}}] at (0,0.25) (a5) {};
        \node[roundnode, fill=black, label=0:{\footnotesize{$u_6$}}] at (0.75,-0.75) (a6) {};
        \node[roundnode, fill=black, label=270:{\footnotesize{$u_7$}}] at (0,-1.75) (a7) {};
        
            \draw[-] (aa0) to node[midway, below]{\footnotesize{$1$}} (aa1);
            \draw[-] (aa0) to node[midway, right]{\footnotesize{$4$}} (a1);
            \draw[-] (aa0) to node[midway, below]{\footnotesize{$9$}} (a2);
            \draw[-] (aa0) to node[midway, right]{\footnotesize{$2$}} (a3);
            \draw[-] (aa1) to node[midway, below]{\footnotesize{$1$}} (aa2);
            \draw[-] (aa1) to node[midway, right]{\footnotesize{$1$}} (a4);
            \draw[-] (aa2) to node[midway, left]{\footnotesize{$3$}} (a5);
            \draw[-] (aa2) to node[midway, above]{\footnotesize{$5$}} (a6);
            \draw[-] (aa2) to node[midway, left]{\footnotesize{$7$}} (a7);
          \end{tikzpicture}
\vspace{-25pt}
\subcaption{}
\end{center}
\end{minipage}
\begin{minipage}{0.35\linewidth}
\begin{center}
\begin{tikzpicture}[roundnode/.style={circle, draw=black, fill=black, inner sep=0pt, minimum size=1.5mm}]

        \node[roundnode, label=180:{\footnotesize{$u_1$}}] at (120:1.05cm) (a1) {};
        \node[roundnode, label=0:{\footnotesize{$u_2$}}] at (70:1.05cm) (a2) {};
        \node[roundnode, label=0:{\footnotesize{$u_3$}}] at (20:1.05cm) (a3) {};
        \node[roundnode, label=-0:{\footnotesize{$u_4$}}] at (-30:1.05cm) (a4) {};
        \node[roundnode, label=270:{\footnotesize{$u_5$}}] at (-80:1.05cm) (a5) {};
        \node[roundnode, label=-180:{\footnotesize{$u_6$}}] at (-130:1.05cm) (a6) {};
        \node[roundnode, label=-180:{\footnotesize{$u_7$}}] at (-180:1.05cm) (a7) {};
        
            \draw[-] (a1) -- (a2);
            \draw[-] (a1) -- (a5);
            \draw[-] (a1) -- (a6);
            \draw[-] (a1) -- (a7);
            \draw[-] (a2) -- (a3);
            \draw[-] (a2) -- (a4);
            \draw[-] (a3) -- (a6);
            \draw[-] (a3) -- (a7);
            \draw[-] (a5) -- (a7);
            \draw[-] (a4) -- (a7); 
            \draw[-] (a6) -- (a7);
          \end{tikzpicture}
\vspace{-17pt}
\subcaption{}
\end{center}
\end{minipage}
\caption{(a) An edge-weighted tree $T$, (b) a PCG $G$ of $T$ for $d_{min} = 9$, $d_{max} = 13$}
\label{fig:fig1}
\end{figure}
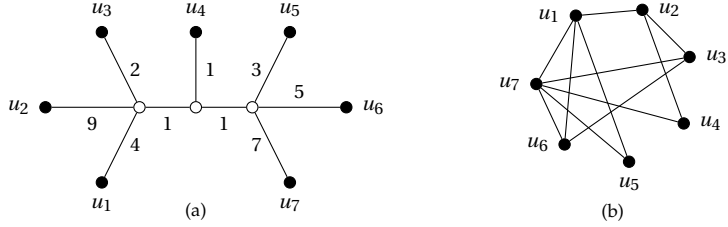

Let $T$ be an edge-weighted tree of non-negative real-valued edge weights with leaf set $L$. Let $d_T(u, v)$ be the sum of the weights of the edges on the path from $u$ to $v$ in $T$, and $d_{min}$ and $ d_{max}$ be two non-negative real numbers where $d_{min} \leq d_{max}$. The \textit{pairwise compatibility graph} (PCG) of $T$ is a graph $G=(V,E)$, where each leaf $u \in L$ corresponds to a vertex $u \in V$ and there is an edge $(u, v) \in E$ if and only if $d_T(u, v)$ lies within the interval $[d_{min}, d_{max}]$. $T$ is called the \textit{pairwise compatibility tree} (PCT) or the \textit{witness tree} of $G$. Figure \ref{fig:fig1}(a) represents an edge-weighted tree $T$ and Figure \ref{fig:fig1}(b) represents a pairwise compatibility graph $G$ of $T$ for $d_{min} = 9$ and $d_{max} = 13$. 

Kearney et al. \cite{kearney2003efficient} introduced the concept of PCGs while working on the phylogenetic tree reconstruction problem. They also found a relation about the clique problem and the PCT construction problem for a given PCG. In a recent research work, Long et al. \cite{long2020exact} showed that PCGs can be applied to describe rare evolutionary events and scenarios with horizontal gene transfer. 

Kearney et al. \cite{kearney2003efficient} conjectured that all undirected graphs are PCGs. But Yanhaona et al. \cite{doi:10.1142/S1793830910000917} refuted the conjecture by constructing a fifteen vertex bipartite graph which is not a PCG. The computational complexity of recognizing a graph as a PCG is not known till now. Hossain et al. \cite{JGAA419} gave a necessary condition and a sufficient condition for a graph to be PCG. They also identified four graph classes that are intermediate to the necessary condition and the sufficient condition and gave examples of two graphs that belong to two of these classes but are not PCGs. But finding examples of graphs that belong to the remaining two classes and are not PCGs is still an open problem. However, any complete characterization of PCGs is still unknown. Xiao et al. \cite{xiao2018characterizing} gave a complete characterization for a graph to be a PCG that admits a star as its witness tree. Calamoneri et al. \cite{CalamoneriCaterpillar} found some characteristics of a PCG that admits a caterpillar as its witness tree. In another work, Calamoneri et al. \cite{CALAMONERI201323} analyzed the closure properties of PCGs under different graph operations. 

Researchers have also worked on some relaxations on the constraints of PCGs. Their works have lead to some subclasses of PCG \cite{CALAMONERI201323, NISHIMURA200269, brandstadt1999distance, brandstadt2008structure, brandstadt2008ptolemaic, BRANDSTADT2010897, kennedy2006strictly, lin2000phylogenetic, nevries2016towards, fellows2008leaf, calamoneri2014relating}. {\color{black} Ahmed et al. \cite{ahmed2017multi} defined multi-interval PCGs as a superclass of PCG. A graph $G$ is a $k$-interval PCG of an edge weighted tree $T$ for mutually exclusive intervals $I_1, I_2, \dots , I_k$ of non-negative real numbers when each vertex of $G$ corresponds to a leaf of $T$ and there is an edge between two vertices in $G$ if the distance between their corresponding leaves lies in $I_1 \cup I_2 \cup \dots \cup I_k$. }

Constructing a PCT from a given graph $G$ for a suitable $d_{min}$ and $d_{max}$ is an interesting problem and a lot of research works have been conducted on identifying different classes of graphs as PCGs {\color{black}like complete graphs, trees, cycles, block-cycle graphs, interval graphs, threshold graphs, threshold tolerance graphs, split matching graphs, ladder graphs, outer subdivision of ladder graphs,  triangle-free outer planar 3-graphs, Dilworth-1 graphs, Dilworth-2 graphs, every graph of at most seven vertices etc} \cite{doi:10.1142/S1793830910000917, CalamoneriCaterpillar, CALAMONERI201323, salma2012triangle, 10.1093/comjnl/bxs087, CALAMONERI201434, yanhaona2009pairwise, rahmanSurvey2020, calamoneri2016pairwise}.  Subsequently, some other graphs are also discovered which are not PCGs like {\color{black}  wheel graphs with at least 9 vertices, certain bipartite graphs etc} \cite{DUROCHER201578, calamoneri2013pairwise, baiocchi2018graphs}. {\color{black} Azam et al. \cite{azamNumber} formulated a method to enumerate all PCGs for a given number of vertices and found that all but seven graphs of eight vertices are PCGs. Later, Azam et al. \cite{pqw} gave an algorithm to enumerate all minimal nonPCGs for a given number of vertices, where a minimal nonPCG is such a graph that is not a PCG but each of its induced subgraph is.} 

It is still unknown whether grid graphs are PCGs or not. But it has been proven that disk graphs, a superclass of grid graphs, are not PCGs \cite{rahmanSurvey2020} and ladder graphs, a subclass of grid graphs, are PCGs \cite{salma2012triangle}. {\color{black} The PCT of a ladder graph, as constrcuted by \cite{salma2012triangle}, is a centepede where the order of leaves connected to the spine is the row major order of the original ladder graph. Using a similar idea, grid graphs have been shown as 2-interval PCGs by Papan et al. \cite{papan}. In this paper, 
we take a different route to prove that grid graphs are PCGs.} We construct examples of the remaining two intermediate classes introduced by Hossain et al.\ \cite{JGAA419} and show that these two graphs are not PCGs. 

The remaining of the paper is organized as follows. Section \ref{sec:preliminaries} gives the definitions of the terms that we use afterwards and revisits some previous findings regarding the necessary and the sufficient condition of PCGs. In section \ref{sec:grid_graphs}, we show that grid graphs are PCGs. In section \ref{sec:neg_examples}, we construct novel negative examples for two intermediate graph classes. Finally, section \ref{sec:conclusion} concludes our paper with discussions. 

\section{Preliminaries}
\label{sec:preliminaries}
In this section, we first define some terminologies which will be used throughout this paper, referring to \cite{rahman2017basic} for graph theoretic terminologies which are not defined here. {\color{black} Afterwards, we discuss some previous results relevant to our works.}

\subsection{{\color{black}Definitions}}

Let $G = (V,E)$ be a simple, undirected graph with vertex set $V$ and edge set $E$. 
An \textit{isomorphism} between two graphs $G_1 = (V_1, E_1)$ and $G2 = (V_2, E_2)$ is a one-to-one correspondence between the vertices in $V_1$ and $V_2$ such that the number of edges between any two vertices in $V_1$ is equal to the number of edges between the corresponding two vertices in $V_2$. If there is an isomorphism between two graphs $G_1$ and $G_2$, then we say that $G_1$ is \textit{isomorphic} to $G_2$. The \textit{complement} of a graph $G = (V, E)$ is a graph $G^c = (V, E^c)$ on the same vertex set $V$ such that for any two vertices $u,v \in V$, $(u,v) \in E^c$ if and only if $(u, v) \notin E$. A \textit{caterpillar} is a tree in which the deletion of all leaves produces a path. The path is called the \textit{spine} of the caterpillar. 
A \textit{grid graph}, $G_{k,l}$, 
is such a graph that vertices correspond to the grid points of a $k \times l$ grid in the plane and edges correspond to the grid line segments between consecutive grid points. For $1 \leq x \leq k$, $1 \leq y \leq l$, we denote the vertex corresponding to the grid point $(x,y)$ by $u_{x,y}$. The following fact holds for grid graphs.

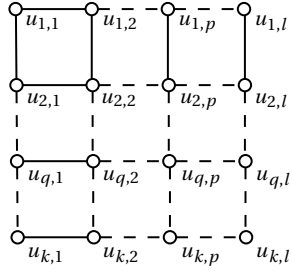
\begin{figure}[H]
    \centering
    \tikzset{every picture/.style={line width=0.75pt}} 

\begin{tikzpicture}[x=0.75pt,y=0.75pt,yscale=-0.6,xscale=0.6]

\draw   (56,45) .. controls (56,42.13) and (58.33,39.8) .. (61.2,39.8) .. controls (64.07,39.8) and (66.4,42.13) .. (66.4,45) .. controls (66.4,47.87) and (64.07,50.2) .. (61.2,50.2) .. controls (58.33,50.2) and (56,47.87) .. (56,45) -- cycle ;
\draw   (120,45) .. controls (120,42.13) and (122.33,39.8) .. (125.2,39.8) .. controls (128.07,39.8) and (130.4,42.13) .. (130.4,45) .. controls (130.4,47.87) and (128.07,50.2) .. (125.2,50.2) .. controls (122.33,50.2) and (120,47.87) .. (120,45) -- cycle ;
\draw   (184,45) .. controls (184,42.13) and (186.33,39.8) .. (189.2,39.8) .. controls (192.07,39.8) and (194.4,42.13) .. (194.4,45) .. controls (194.4,47.87) and (192.07,50.2) .. (189.2,50.2) .. controls (186.33,50.2) and (184,47.87) .. (184,45) -- cycle ;
\draw   (248,45) .. controls (248,42.13) and (250.33,39.8) .. (253.2,39.8) .. controls (256.07,39.8) and (258.4,42.13) .. (258.4,45) .. controls (258.4,47.87) and (256.07,50.2) .. (253.2,50.2) .. controls (250.33,50.2) and (248,47.87) .. (248,45) -- cycle ;
\draw    (66.4,45) -- (120,45) ;
\draw  [dash pattern={on 4.5pt off 4.5pt}]  (130.4,45) -- (184,45) ;
\draw  [dash pattern={on 4.5pt off 4.5pt}]  (194.4,45) -- (248,45) ;
\draw   (61.91,103.8) .. controls (64.79,103.78) and (67.13,106.09) .. (67.15,108.96) .. controls (67.18,111.83) and (64.87,114.18) .. (61.99,114.2) .. controls (59.12,114.22) and (56.78,111.91) .. (56.75,109.04) .. controls (56.73,106.17) and (59.04,103.82) .. (61.91,103.8) -- cycle ;
\draw   (62.41,167.8) .. controls (65.28,167.78) and (67.62,170.09) .. (67.65,172.96) .. controls (67.67,175.83) and (65.36,178.18) .. (62.49,178.2) .. controls (59.61,178.22) and (57.27,175.91) .. (57.25,173.04) .. controls (57.22,170.17) and (59.53,167.82) .. (62.41,167.8) -- cycle ;
\draw   (62.9,231.8) .. controls (65.77,231.78) and (68.12,234.09) .. (68.14,236.96) .. controls (68.16,239.83) and (65.85,242.17) .. (62.98,242.2) .. controls (60.11,242.22) and (57.76,239.91) .. (57.74,237.04) .. controls (57.72,234.17) and (60.03,231.82) .. (62.9,231.8) -- cycle ;
\draw    (61.5,50.2) -- (61.91,103.8) ;
\draw  [dash pattern={on 4.5pt off 4.5pt}]  (61.99,114.2) -- (62.41,167.8) ;
\draw  [dash pattern={on 4.5pt off 4.5pt}]  (62.49,178.2) -- (62.9,231.8) ;
\draw    (67.15,108.96) -- (120.75,108.96) ;
\draw   (120.75,108.96) .. controls (120.75,106.09) and (123.08,103.76) .. (125.95,103.76) .. controls (128.83,103.76) and (131.15,106.09) .. (131.15,108.96) .. controls (131.15,111.83) and (128.83,114.16) .. (125.95,114.16) .. controls (123.08,114.16) and (120.75,111.83) .. (120.75,108.96) -- cycle ;
\draw  [dash pattern={on 4.5pt off 4.5pt}]  (131.15,108.96) -- (184.75,108.96) ;
\draw   (184.75,108.96) .. controls (184.75,106.09) and (187.08,103.76) .. (189.95,103.76) .. controls (192.83,103.76) and (195.15,106.09) .. (195.15,108.96) .. controls (195.15,111.83) and (192.83,114.16) .. (189.95,114.16) .. controls (187.08,114.16) and (184.75,111.83) .. (184.75,108.96) -- cycle ;
\draw  [dash pattern={on 4.5pt off 4.5pt}]  (195.15,108.96) -- (248.75,108.96) ;
\draw   (248.75,108.96) .. controls (248.75,106.09) and (251.08,103.76) .. (253.95,103.76) .. controls (256.83,103.76) and (259.15,106.09) .. (259.15,108.96) .. controls (259.15,111.83) and (256.83,114.16) .. (253.95,114.16) .. controls (251.08,114.16) and (248.75,111.83) .. (248.75,108.96) -- cycle ;
\draw    (67.65,172.96) -- (121.25,172.96) ;
\draw   (121.25,172.96) .. controls (121.25,170.09) and (123.57,167.76) .. (126.45,167.76) .. controls (129.32,167.76) and (131.65,170.09) .. (131.65,172.96) .. controls (131.65,175.83) and (129.32,178.16) .. (126.45,178.16) .. controls (123.57,178.16) and (121.25,175.83) .. (121.25,172.96) -- cycle ;
\draw   (185.25,172.96) .. controls (185.25,170.09) and (187.57,167.76) .. (190.45,167.76) .. controls (193.32,167.76) and (195.65,170.09) .. (195.65,172.96) .. controls (195.65,175.83) and (193.32,178.16) .. (190.45,178.16) .. controls (187.57,178.16) and (185.25,175.83) .. (185.25,172.96) -- cycle ;
\draw  [dash pattern={on 4.5pt off 4.5pt}]  (131.65,172.96) -- (185.25,172.96) ;
\draw  [dash pattern={on 4.5pt off 4.5pt}]  (195.65,172.96) -- (249.25,172.96) ;
\draw   (249.25,172.96) .. controls (249.25,170.09) and (251.57,167.76) .. (254.45,167.76) .. controls (257.32,167.76) and (259.65,170.09) .. (259.65,172.96) .. controls (259.65,175.83) and (257.32,178.16) .. (254.45,178.16) .. controls (251.57,178.16) and (249.25,175.83) .. (249.25,172.96) -- cycle ;
\draw    (68.14,236.96) -- (121.74,236.96) ;
\draw   (121.74,236.96) .. controls (121.74,234.09) and (124.07,231.76) .. (126.94,231.76) .. controls (129.81,231.76) and (132.14,234.09) .. (132.14,236.96) .. controls (132.14,239.83) and (129.81,242.16) .. (126.94,242.16) .. controls (124.07,242.16) and (121.74,239.83) .. (121.74,236.96) -- cycle ;
\draw  [dash pattern={on 4.5pt off 4.5pt}]  (132.14,236.96) -- (185.74,236.96) ;
\draw  [dash pattern={on 4.5pt off 4.5pt}]  (196.14,236.96) -- (249.74,236.96) ;
\draw   (185.74,236.96) .. controls (185.74,234.09) and (188.07,231.76) .. (190.94,231.76) .. controls (193.81,231.76) and (196.14,234.09) .. (196.14,236.96) .. controls (196.14,239.83) and (193.81,242.16) .. (190.94,242.16) .. controls (188.07,242.16) and (185.74,239.83) .. (185.74,236.96) -- cycle ;
\draw   (249.74,236.96) .. controls (249.74,234.09) and (252.07,231.76) .. (254.94,231.76) .. controls (257.81,231.76) and (260.14,234.09) .. (260.14,236.96) .. controls (260.14,239.83) and (257.81,242.16) .. (254.94,242.16) .. controls (252.07,242.16) and (249.74,239.83) .. (249.74,236.96) -- cycle ;
\draw    (125.2,50.2) -- (125.2,103.76) ;
\draw  [dash pattern={on 4.5pt off 4.5pt}]  (125.95,114.16) -- (125.95,167.76) ;
\draw  [dash pattern={on 4.5pt off 4.5pt}]  (126.45,178.16) -- (126.45,231.76) ;
\draw    (189.2,50.2) -- (189.2,103.76) ;
\draw  [dash pattern={on 4.5pt off 4.5pt}]  (189.95,114.16) -- (189.95,167.72) ;
\draw  [dash pattern={on 4.5pt off 4.5pt}]  (190.94,178.2) -- (190.94,231.76) ;
\draw    (253.95,50.2) -- (253.95,103.76) ;
\draw  [dash pattern={on 4.5pt off 4.5pt}]  (253.95,114.16) -- (253.95,167.72) ;
\draw  [dash pattern={on 4.5pt off 4.5pt}]  (254.45,178.16) -- (254.45,231.72) ;

\draw (66.2,48) node [anchor=north west][inner sep=0.75pt]   [align=left] {\footnotesize{$u{}_{1,1}$}};
\draw (131.2,48) node [anchor=north west][inner sep=0.75pt]   [align=left] {\footnotesize{$u{}_{1,2}$}};
\draw (195.2,48) node [anchor=north west][inner sep=0.75pt]   [align=left] {\footnotesize{$u{}_{1,p}$}};
\draw (259.2,48) node [anchor=north west][inner sep=0.75pt]   [align=left] {\footnotesize{$u{}_{1,l}$}};
\draw (67.2,113) node [anchor=north west][inner sep=0.75pt]   [align=left] {\footnotesize{$u{}_{2,1}$}};
\draw (130.2,113) node [anchor=north west][inner sep=0.75pt]   [align=left] {\footnotesize{$u{}_{2,2}$}};
\draw (196.2,113) node [anchor=north west][inner sep=0.75pt]   [align=left] {\footnotesize{$u{}_{2,p}$}};
\draw (259.2,113) node [anchor=north west][inner sep=0.75pt]   [align=left] {\footnotesize{$u{}_{2,l}$}};
\draw (68.2,178) node [anchor=north west][inner sep=0.75pt]   [align=left] {\footnotesize{$u{}_{q,1}$}};
\draw (131.2,178) node [anchor=north west][inner sep=0.75pt]   [align=left] {\footnotesize{$u{}_{q,2}$}};
\draw (197.2,178) node [anchor=north west][inner sep=0.75pt]   [align=left] {\footnotesize{$u{}_{q,p}$}};
\draw (260.2,178) node [anchor=north west][inner sep=0.75pt]   [align=left] {\footnotesize{$u{}_{q,l}$}};
\draw (68.2,243) node [anchor=north west][inner sep=0.75pt]   [align=left] {\footnotesize{$u{}_{k,1}$}};
\draw (131.2,243) node [anchor=north west][inner sep=0.75pt]   [align=left] {\footnotesize{$u{}_{k,2}$}};
\draw (197.2,243) node [anchor=north west][inner sep=0.75pt]   [align=left] {\footnotesize{$u{}_{k,p}$}};
\draw (260.2,243) node [anchor=north west][inner sep=0.75pt]   [align=left] {\footnotesize{$u{}_{k,l}$}};

\end{tikzpicture}
    \vspace{-20pt}
    \caption{A grid graph $G_{k,l}$}
    \label{grid_graph}
\end{figure}

\begin{fact}
\label{fact:grid_def}
In a grid graph $G_{k,l}$, there is an edge $(u_{x,y}, u_{x',y'})$ if and only if either of the following two conditions is satisfied\\
(C1): $x=x'$ and $|y-y'|=1$  or \\
(C2): $|x-x'|=1$ and $y=y'$.
\end{fact}


\subsection{{\color{black}Some Previous Findings}}

Hossain et al. \cite{JGAA419} gave a necessary and a sufficient condition for a graph to be a PCG.

\begin{lemma}[A Necessary Condition for PCGs]
Let $G$ be a graph. Let $H_1$ and $H_2$ be two disjoint induced subgraphs of $G^c$. If each of $H_1$ and $H_2$ is either a chordless cycle of at least four vertices
or $C_n^c$ for $n \geq 5$, then $G$ is not a PCG. 
\label{theor_nec}
\end{lemma}

\begin{lemma}[A Sufficient Condition for PCGs]
Let $G$ be a graph. If $G^c$ has no cycle, then $G$ is a PCG.

\label{theor_suf}
\end{lemma}


In an attempt to characterize the intermediate graphs, Hossain et al. \cite{JGAA419} proposed four classes. 

\begin{description}
  \item [$\mathcal{G}^1$] A graph $G$ belongs to $\mathcal{G}^1$ if $G^c$ does not contain any chordless cycle.
  \item [$\mathcal{G}^2$] A graph $G$ belongs to $\mathcal{G}^2$ if $G^c$ consists of two induced chordless cycles, where the cycles share some common vertices.
  \item [$\mathcal{G}^3$] A graph $G$ belongs to $\mathcal{G}^3$ if $G^c$ consists of two induced chordless cycles and some edges that are incident to both cycles.

  \item [$\mathcal{G}^4$] A graph $G$ belongs to $\mathcal{G}^4$ if $G^c$ contains only one induced chordless cycle.
\end{description}

They provided examples of graphs in $\mathcal{G}^1$ and $\mathcal{G}^3$ classes that are not PCGs.
Whether all graphs contained in $\mathcal{G}^2$ and $\mathcal{G}^4$ classes are PCGs remained an open problem \cite{JGAA419,rahmanSurvey2020}. We construct examples in both $\mathcal{G}^2$ and $\mathcal{G}^4$ and prove them not to be PCGs. To facilitate our construction, we need to revisit the negative example given by Yanhaona et al \cite{doi:10.1142/S1793830910000917}. They showed a bipartite graph $H$ of 15 vertices that is not a PCG. The vertices of the graph $H$ can be divided into two disjoint independent sets $A = \{a_1, a_2, ..., a_5\}$ and $B = \{b_1, b_2, ..., b_{10}\}$. We call the sets $A$ and $B$ the independent sets of $H$. Each vertex in $A$ has edges with six distinct vertices in $B$ and each vertex in $B$ has edges with three distinct vertices in $A$. The neighborhood of each vertex is unique. 
\begin{figure}
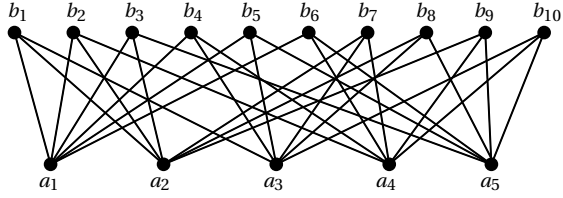

    \centering
    \include{first_graph_non_PCG}
    \vspace{-20pt}
    \caption{The graph $H$, the first graph that was proved not to be a PCG.}
    \label{H}
\end{figure}

Hossain et al. \cite{JGAA419} revisited the proof of $H$ not being a PCG and reached the following conclusion. 

\begin{lemma}
  
Let $\mathcal{H}$ be the class of graphs consisting of all the graphs obtained from $H$ by adding edges between the vertices in the same independent set of $H$. Then none of the graphs in $\mathcal{H}$ is a PCG.

\label{fl1}
\end{lemma}

We use Lemma \ref{fl1} for constructing our negative example in Section \ref{subsection:neg_example_1}.


\section{Grid Graphs are PCGs}
\label{sec:grid_graphs}

In this section we prove that grid graphs are PCGs by providing an algorithm to construct a PCT of a given grid graph $G_{k,l}$.

\begin{figure}[H]
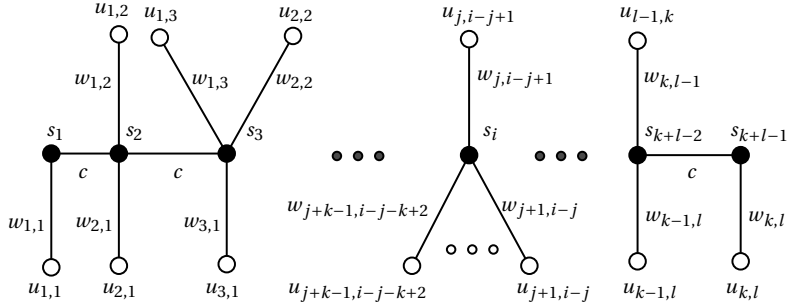

    \centering
    \include{PCT}
    \vspace{-20pt}
    \caption{A PCT $T_g$ of $G_{k,l}$}
    \label{pct_fig}
\end{figure}

Figure \ref{grid_graph} shows a grid graph $G_{k,l}$. If $k=2$ or $l=2$, $G_{k,l}$ becomes a ladder graph which is already proven to be a PCG. As $G_{k,l}$ is isomorphic to $G_{l,k}$, here we can consider $3 \leq k \leq l$. We define a set of vertices, $D_i$, as a \textit{diagonal set} of $G_{k,l}$ such that $D_i$ consists of each vertex $u_{x,y} \in V(G_{k,l})$ where $i=x+y-1$. It is trivial that the total number of diagonal sets in $G_{k,l}$ is $k+l-1$ and so $1 \leq i \leq k+l-1$. As $k \leq l$, it is obvious that the total number of vertices in a diagonal set of $G_{k,l}$ is at most $k$. One can observe that the diagonal sets induce a partition of the vertex set of $G_{k,l}$. That is, $\bigcup_{i=1}^{k+l-1}D_i=V(G_{k,l})$, and for $1 \leq i,j \leq k+l-1$ and $i \neq j$, $D_i \cap D_j = \emptyset$.  In addition, $|y-x| < l$ for each $1 \leq x \leq k$ and $1 \leq y \leq l$.

We now give an algorithm to construct an edge-weighted caterpillar $T_g$ corresponding to a grid graph $G_{k,l}$. For each $1 \leq i \leq k+l-1$, we set a vertex $s_i$ on the spine of $T_g$ corresponding to diagonal set $D_i$ in $G_{k,l}$. Then for each vertex $u_{x,y}$ of $G_{k,l}$ in diagonal set $D_i$, we draw a leaf vertex $u_{x,y}$ and an edge $(s_i, u_{x,y})$ in $T_g$. Figure \ref{pct_fig} shows the construction of $T_g$. We denote the weight of the edge $(s_i, u_{x,y}) \in E(T_g)$ by $w_{x, y}$ where $w_{x,y} = r-(-1)^i(y-x)$ for each $1 \leq i \leq k+l-1$, $1 \leq x \leq k$, $1 \leq y \leq l$ and $r = k+l+1$. Since $|y-x| < l < r$, all values of $w_{x,y}$ are non-negative. For each $1 \leq i \leq k+l-2$, we set the weight of the edge $(s_i, s_{i+1}) \in E(T_g)$ as $c$ where $c = 4k+4l+4$. We now prove that $T_g$ is a PCT of $G_{k,l}$ as in the following theorem.

\begin{theorem}
    $G_{k,l}$ is a PCG of $T_g$ for $d_{min} = 2r+c-1$ and $d_{max} = 2r+c+1$ where $r = k+l+1$ and $c = 4k+4l+4$. 
\end{theorem}

\begin{proof}
Since the diagonal sets induce a partition of the vertex set of $G_{k,l}$, our construction implies that there is a bijection between the set of leaves of $T_g$ and the set of vertices of $G_{k,l}$. Let $u_{x,y}$ and $u_{x',y'}$ be two leaves in $T_g$.

It is sufficient to prove the following two claims to establish that $T_g$ is a PCT of $G_{k,l}$ with the mentioned values of $d_{min}$ and $d_{max}$.

\begin{description}
    \item [Claim 1.] If $u_{x,y}$ and $u_{x',y'}$ are not adjacent in $G_{k,l}$, then \\$d_{T_g}(u_{x,y},u_{x',y'}) \notin [d_{min}, d_{max}]$.
    \item [Claim 2.] If $u_{x,y}$ and $u_{x',y'}$ are adjacent in $G_{k,l}$, then \\$d_{T_g}(u_{x,y},u_{x',y'}) \in [d_{min}, d_{max}]$.
\end{description}

We first prove Claim 1.

 When the leaves $u_{x,y}$ and $u_{x',y'}$ belong to the same diagonal set $D_i$, the corresponding vertices are not adjacent in $G_{k,l}$. Our construction entails that $x + y = x' + y' = i+1$ and  $d_{T_g}(u_{x,y},u_{x',y'}) =  w_{x,y} + w_{{x'},{y'}} = r - (-1)^i(y-x) + r - (-1)^i(y'-x')  \leq 2r + |y-x| + |y'-x'| < 2r+l+l < 2r+4k+4l+3 = 2r+c-1 = d_{min}$. Since $d_{min} < d_{max}$,  $d_{T_g}(u_{x,y},u_{x',y'}) \notin [d_{min}, d_{max}]$.

Now let the leaves $u_{x,y}$ and $u_{x',y'}$ belong to the different diagonal sets $D_i$ and $D_{i'}$ respectively. Let $s_i$ and $s_{i'}$ be the nodes on the spine of $T_g$ which are adjacent to $u_{x,y}$ and $u_{x',y'}$, respectively. Then, from the construction of caterpillar $T_g$, for each $1 \leq x,x' \leq k$, $1 \leq y,y' \leq l$, $1\leq i,i' \leq k+l-1$ and $i \neq i'$, we can write,
\begin{equation}
\label{eq1}
d_{T_g}(u_{x,y},u_{x',y'}) =  w_{x,y} + w_{s_i, s_{i'}} + w_{{x'},{y'}}.  
\end{equation}

Since $w_{x,y} = r-(-1)^i(y-x)$, $w_{s_i, s_{i'}} = |i-i'|c$ and $w_{{x'},{y'}} = r-(-1)^{i'}(y'-x')$, equation (\ref{eq1}) can be written as:
\begin{equation}
\label{eq2}
d_{T_g}(u_{x,y},u_{{x'},{y'}}) = r-(-1)^i(y-x) +  |i-i'|c + r-(-1)^{i'}({y'}-{x'}).
\end{equation}

Let $|i-i'|\geq2$. By Fact \ref{fact:grid_def}, the leaves $u_{x,y}$ and $u_{x',y'}$ are not adjacent. Then \\$d_{T_g}(u_{x,y},u_{{x'},{y'}}) \geq r - |y-x| - |y'-x'| + 2c > r - l - l + c + c = k+l+1 -2l + 4k+4l+4 + c = 5k+3l+c+5 > 2k+2l+2+c+1 = 2r+c+1 = d_{max}$. Since $d_{min} < d_{max}$,  $d_{T_g}(u_{x,y},u_{x',y'}) \notin [d_{min}, d_{max}]$.

What remains for proving Claim 1 is the case when $u_{x,y}$ and $u_{x',y'}$ are not adjacent in $G_{k, l}$ and $|i-i'| = 1$. This part will be proved by contraposition.  Here $i$ and $i'$ has opposite parity. Without the loss of generality, let $i$ be even and $i'$ be odd. It follows that $d_{T_g}(u_{x,y},u_{{x'},{y'}}) = r-(y-x) + c + r+({y'}-{x'}) = 2r+c+(y'-y)+(x-x')$. Here $d_{min} = 2r+c-1$ and $d_{max} = 2r+c+1$. Let $(y'-y)+(x-x')=0$. Then $|i-i'|=|(x-x')-(y'-y)|=|-2(y'-y)| = 2 |y'-y| =1$ which is not realizable. Again, for $u_{x,y}$ and $u_{x',y'}$ to be adjacent in $G_{k, l}$, either (C1) or (C2) from Fact \ref{fact:grid_def} must hold.

Thus, to prove Claim 1, it suffices to prove that if $(y'-y)+(x-x') = \pm 1$, either (C1) or (C2) holds true. We will consider each value of $(y'-y)+(x-x')$ for \begin{equation}
        \label{app1}
        i-i' = (x+y) - (x'+y') = -(y'-y)+(x-x') =1
\end{equation}

Let $(y'-y)+(x-x') = -1$. Adding with Equation \ref{app1}, we get $2(x-x')=0 \implies x=x'$ and $y'-y = -1$ and (C1) is satisfied. Now let $(y'-y)+(x-x') = 1$. Adding with Equation \ref{app1}, we get $2(x-x')=2 \implies x-x'=1$ and $y'=y$ and (C2) is satisfied.\\

Now we will consider similar cases for 
\begin{equation}
        \label{app2}
        i-i' = (x+y) - (x'+y') = -(y'-y)+(x-x') = -1
\end{equation}

Let $(y'-y)+(x-x') = -1$. Adding with Equation \ref{app2}, we get $2(x-x')=-2 \implies x-x'=-1$ and $y'=y$ and (C2) is satisfied. Now let $(y'-y)+(x-x') = 1$. Adding with Equation \ref{app2}, we get $2(x-x')=0 \implies x=x'$ and $y'-y=1$ and (C1) is satisfied. Thus Claim 1 is proved.



We now prove Claim 2. Here we assume that $u_{x,y}$ and $u_{x',y'}$ are adjacent in $G_{k, l}$. Then either (C1) or (C2) from Fact 1 holds. In either case $|i-i'|=1$. Without the loss of generality, let us again assume $i$ is even and $i'$ is odd. From Equation \ref{eq2}, $$d_{T_g}(u_{x,y},u_{{x'},{y'}}) = r-(y-x) +  c + r+({y'}-{x'}) = 2r + c + (x-x') + (y'-y)$$

Since either (C1) or (C2) holds, $d_{T_g}(u_{x,y},u_{{x'},{y'}}) = 2r+c \pm 1 \in [d_{min}, d_{max}]$.
\qed \end{proof}

We can construct a PCT of a given grid graph using the above algorithm which clearly runs in linear time. This construction of caterpillar $T_g$ also supports a finding of Calamoneri et al.  \cite{CalamoneriCaterpillar} as $G_{k,l}$ is a triangle-free graph.

\section{Negative Examples}
\label{sec:neg_examples}
In this section, addressing the open problems in \cite{JGAA419}, we show negative examples in some classes of graphs that are intermediate to the classes  defined by the necessary and the sufficient condition given by Hossain et al. \cite{JGAA419}.

\subsection{Negative Example for $\mathcal{G}^2$}
\label{subsection:neg_example_1}

\begin{figure}
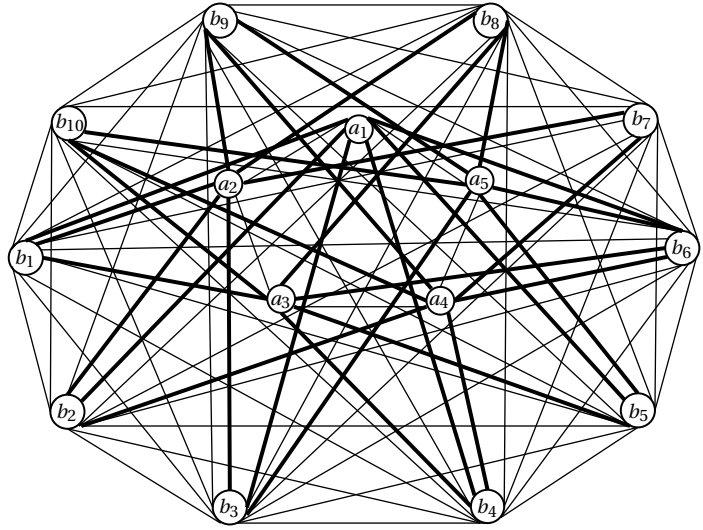

    \centering
    \include{counter_example_G2}
    \vspace{-15pt}
    \caption{A negative example graph $H_2$ in the class $\mathcal{G}^2$. The thick edges were already in $H$. The thin edges are newly introduced.}
    \label{H2}
\end{figure}

We take the graph $H$ (in Figure \ref{H}) which was proved not to be a PCG. We add edges so that the vertices in $A$ form a $5$-cycle and vertices in $B$ form a $10$-clique. Let this graph, presented in Figure \ref{H2}, be $H_2$. By Lemma \ref{fl1}, the graph $H_2$ is not a PCG. We prove the following theorem by showing that the graph $H_2$ is in the class $\mathcal{G}^2$. 

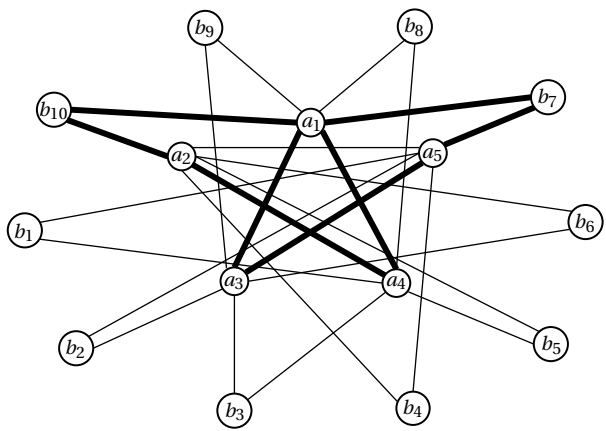
\begin{figure}
    \centering
    \tikzset{every picture/.style={line width=0.75pt}} 

\begin{tikzpicture}[x=0.75pt,y=0.75pt,yscale=-0.7,xscale=0.7]

\draw  [fill={rgb, 255:red, 255; green, 255; blue, 255 }  ,fill opacity=1 ] (329.3,112.6) .. controls (329.3,107.08) and (333.78,102.6) .. (339.3,102.6) .. controls (344.82,102.6) and (349.3,107.08) .. (349.3,112.6) .. controls (349.3,118.12) and (344.82,122.6) .. (339.3,122.6) .. controls (333.78,122.6) and (329.3,118.12) .. (329.3,112.6) -- cycle ;
\draw   (142,103.3) .. controls (142,96.51) and (147.51,91) .. (154.3,91) .. controls (161.09,91) and (166.6,96.51) .. (166.6,103.3) .. controls (166.6,110.09) and (161.09,115.6) .. (154.3,115.6) .. controls (147.51,115.6) and (142,110.09) .. (142,103.3) -- cycle ;
\draw   (121,190.3) .. controls (121,183.51) and (126.51,178) .. (133.3,178) .. controls (140.09,178) and (145.6,183.51) .. (145.6,190.3) .. controls (145.6,197.09) and (140.09,202.6) .. (133.3,202.6) .. controls (126.51,202.6) and (121,197.09) .. (121,190.3) -- cycle ;
\draw   (525,184.3) .. controls (525,177.51) and (530.51,172) .. (537.3,172) .. controls (544.09,172) and (549.6,177.51) .. (549.6,184.3) .. controls (549.6,191.09) and (544.09,196.6) .. (537.3,196.6) .. controls (530.51,196.6) and (525,191.09) .. (525,184.3) -- cycle ;
\draw   (498,94.3) .. controls (498,87.51) and (503.51,82) .. (510.3,82) .. controls (517.09,82) and (522.6,87.51) .. (522.6,94.3) .. controls (522.6,101.09) and (517.09,106.6) .. (510.3,106.6) .. controls (503.51,106.6) and (498,101.09) .. (498,94.3) -- cycle ;
\draw   (251,44.3) .. controls (251,37.51) and (256.51,32) .. (263.3,32) .. controls (270.09,32) and (275.6,37.51) .. (275.6,44.3) .. controls (275.6,51.09) and (270.09,56.6) .. (263.3,56.6) .. controls (256.51,56.6) and (251,51.09) .. (251,44.3) -- cycle ;
\draw   (402,43.3) .. controls (402,36.51) and (407.51,31) .. (414.3,31) .. controls (421.09,31) and (426.6,36.51) .. (426.6,43.3) .. controls (426.6,50.09) and (421.09,55.6) .. (414.3,55.6) .. controls (407.51,55.6) and (402,50.09) .. (402,43.3) -- cycle ;
\draw   (158,275.3) .. controls (158,268.51) and (163.51,263) .. (170.3,263) .. controls (177.09,263) and (182.6,268.51) .. (182.6,275.3) .. controls (182.6,282.09) and (177.09,287.6) .. (170.3,287.6) .. controls (163.51,287.6) and (158,282.09) .. (158,275.3) -- cycle ;
\draw   (272,320.3) .. controls (272,313.51) and (277.51,308) .. (284.3,308) .. controls (291.09,308) and (296.6,313.51) .. (296.6,320.3) .. controls (296.6,327.09) and (291.09,332.6) .. (284.3,332.6) .. controls (277.51,332.6) and (272,327.09) .. (272,320.3) -- cycle ;
\draw   (401,318.6) .. controls (401,311.97) and (406.37,306.6) .. (413,306.6) .. controls (419.63,306.6) and (425,311.97) .. (425,318.6) .. controls (425,325.23) and (419.63,330.6) .. (413,330.6) .. controls (406.37,330.6) and (401,325.23) .. (401,318.6) -- cycle ;
\draw   (500,272.3) .. controls (500,265.51) and (505.51,260) .. (512.3,260) .. controls (519.09,260) and (524.6,265.51) .. (524.6,272.3) .. controls (524.6,279.09) and (519.09,284.6) .. (512.3,284.6) .. controls (505.51,284.6) and (500,279.09) .. (500,272.3) -- cycle ;
\draw  [fill={rgb, 255:red, 255; green, 255; blue, 255 }  ,fill opacity=1 ] (236.3,137) .. controls (236.3,131.48) and (240.78,127) .. (246.3,127) .. controls (251.82,127) and (256.3,131.48) .. (256.3,137) .. controls (256.3,142.52) and (251.82,147) .. (246.3,147) .. controls (240.78,147) and (236.3,142.52) .. (236.3,137) -- cycle ;
\draw  [fill={rgb, 255:red, 255; green, 255; blue, 255 }  ,fill opacity=1 ] (274.3,227) .. controls (274.3,221.48) and (278.78,217) .. (284.3,217) .. controls (289.82,217) and (294.3,221.48) .. (294.3,227) .. controls (294.3,232.52) and (289.82,237) .. (284.3,237) .. controls (278.78,237) and (274.3,232.52) .. (274.3,227) -- cycle ;
\draw  [fill={rgb, 255:red, 255; green, 255; blue, 255 }  ,fill opacity=1 ] (391,228.3) .. controls (391,222.78) and (395.48,218.3) .. (401,218.3) .. controls (406.52,218.3) and (411,222.78) .. (411,228.3) .. controls (411,233.82) and (406.52,238.3) .. (401,238.3) .. controls (395.48,238.3) and (391,233.82) .. (391,228.3) -- cycle ;
\draw  [fill={rgb, 255:red, 255; green, 255; blue, 255 }  ,fill opacity=1 ] (417.3,134.6) .. controls (417.3,129.08) and (421.78,124.6) .. (427.3,124.6) .. controls (432.82,124.6) and (437.3,129.08) .. (437.3,134.6) .. controls (437.3,140.12) and (432.82,144.6) .. (427.3,144.6) .. controls (421.78,144.6) and (417.3,140.12) .. (417.3,134.6) -- cycle ;
\draw [line width=2.25]    (166.6,103.3) -- (329.3,112.6) ;
\draw [line width=2.25]    (284.3,217) -- (331.6,118.6) ;
\draw [line width=2.25]    (347.6,118.6) -- (401,218.3) ;
\draw [line width=0.5]     (345.6,104.6) -- (407.6,52.6) ;
\draw [line width=0.5]     (332.6,104.6) -- (271.6,52.6) ;
\draw [line width=0.5]     (182.6,275.3) -- (275.6,232.6) ;
\draw [line width=0.5]     (284.3,308) -- (284.3,237) ;
\draw [line width=0.5]     (525.6,189.6) -- (294.3,227) ;
\draw [line width=0.5]     (278.6,217.6) -- (263.3,56.6) ;
\draw [line width=2.25]    (292.6,220.6) -- (419.6,141.6) ;
\draw [line width=0.5]     (143.6,184.6) -- (417.3,134.6) ;
\draw [line width=0.5]     (413,306.6) -- (427.3,144.6) ;
\draw [line width=2.25]    (500.6,101.6) -- (435.6,128.6) ;
\draw [line width=0.5]     (417.3,134.6) -- (179.6,266.6) ;
\draw [line width=0.5]     (417.3,130.6) -- (254.6,130.6) ;
\draw [line width=2.25]    (163.6,110.6) -- (236.3,137) ;
\draw [line width=0.5]     (401,218.3) -- (414.3,55.6) ;
\draw [line width=0.5]     (256.3,137) -- (527.6,176.6) ;
\draw [line width=0.5]     (409.6,234.6) -- (500,272.3) ;
\draw [line width=0.5]     (246.3,147) -- (401.6,313.6) ;
\draw [line width=0.5]     (293.6,312.6) -- (393.6,234.6) ;
\draw [line width=0.5]     (256.3,137) -- (504.6,263.6) ;
\draw [line width=0.5]     (143.6,196.6) -- (391,228.3) ;
\draw [line width=2.25]    (253.6,142.6) -- (392.6,222.6) ;
\draw [line width=2.25]    (349.3,112.6) -- (498,94.3) ;

\draw (330,107) node [anchor=north west][inner sep=0.75pt]   [align=left] {\footnotesize{$a_{1}$}};
\draw (236,131) node [anchor=north west][inner sep=0.75pt]   [align=left] {\footnotesize{$a_{2}$}};
\draw (274,221) node [anchor=north west][inner sep=0.75pt]   [align=left] {\footnotesize{$a_{3}$}};
\draw (391,221) node [anchor=north west][inner sep=0.75pt]   [align=left] {\footnotesize{$a_{4}$}};
\draw (417,129) node [anchor=north west][inner sep=0.75pt]   [align=left] {\footnotesize{$a_{5}$}};
\draw (141,95) node [anchor=north west][inner sep=0.75pt]   [align=left] {\footnotesize{$b_{10}$}};
\draw (123,183) node [anchor=north west][inner sep=0.75pt]   [align=left] {\footnotesize{$b_{1}$}};
\draw (527,176) node [anchor=north west][inner sep=0.75pt]   [align=left] {\footnotesize{$b_{6}$}};
\draw (500,86) node [anchor=north west][inner sep=0.75pt]   [align=left] {\footnotesize{$b_{7}$}};
\draw (253,36) node [anchor=north west][inner sep=0.75pt]   [align=left] {\footnotesize{$b_{9}$}};
\draw (404,35) node [anchor=north west][inner sep=0.75pt]   [align=left] {\footnotesize{$b_{8}$}};
\draw (160,267) node [anchor=north west][inner sep=0.75pt]   [align=left] {\footnotesize{$b_{2}$}};
\draw (275,312) node [anchor=north west][inner sep=0.75pt]   [align=left] {\footnotesize{$b_{3}$}};
\draw (403,310) node [anchor=north west][inner sep=0.75pt]   [align=left] {\footnotesize{$b_{4}$}};
\draw (503,264) node [anchor=north west][inner sep=0.75pt]   [align=left] {\footnotesize{$b_{5}$}};

\end{tikzpicture}
    \vspace{-15pt}
    \caption{The complement graph $H_2^c$ of $H_2$. Two chordless cycles $C_1 = \{a_1, a_4, a_2, b_{10}\}$ and $C_2 = \{a_1, a_3, a_5, b_7\}$ share a common vertex $a_1$.}
    \label{H2c}
\end{figure}

\begin{theorem}
    The class $\mathcal{G}^2$ contains some graphs which are not PCGs. 
\end{theorem} 
\begin{proof}
    We already have a graph $H_2$ which is not a PCG. It is sufficient to show $H_2$ is in $\mathcal{G}^2$, i. e., its complement graph $H_2^c$ has two chordless cycles with some common vertices. 
    The complement graph $H_2^c$ is shown in Figure \ref{H2c} where there are two chordless cycles $C_1 = \{a_1, a_4, a_2, b_{10}\}$ and $C_2 = \{a_1, a_3, a_5, b_7\}$. Since two cycles have the vertex $a_1$ in common, the graph $H_4$ belongs to the class $\mathcal{G}^2$. Since the graph $H_4$ is not a PCG, the proof is completed.  
\qed \end{proof}    
    
\subsection{Negative Example for $\mathcal{G}^4$} 

We now take the graph $H_1$ proven not to be a PCG as a negative example of the class $\mathcal{G}^1$ in \cite{rahmanSurvey2020}. Like $H$, it has two sets of vertices $A = \{a_1, a_2, ..., a_5\}$ and $B = \{b_1, b_2, ..., b_{10}\}$. The vertices in $A$ form a 5-clique. We introduce a new set of vertices $C = \{c_1, c_2, ..., c_5\}$ that form a $C_5^c$, the complement graph of a 5-cycle. In addition, from every vertex $v$ in the set $A \cup B$, we add an edge to every vertex in $C$. We call the resultant graph $H_4$. To prove that the graph $H_4$ is not a PCG, we restate the following lemma from \cite{CALAMONERI201323}.  

\begin{figure}[H]
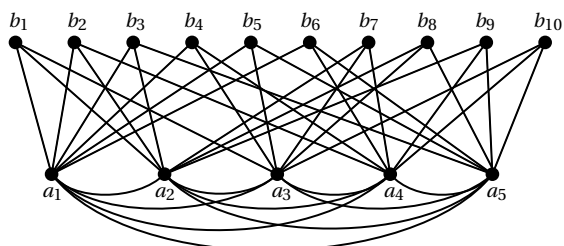

    \centering
    \include{H1}
    \vspace{-30pt}
    \caption{The graph $H_1$ that was proven not to be a PCG.}
    \label{H1}
\end{figure}

\begin{lemma}
     If a graph $G$ is not a PCG, then every
graph containing $G$ as an induced subgraph is also not a PCG. 
\label{fl4}
\end{lemma}
\begin{corollary}
    The graph $H_4$ is not a PCG.
\label{fl4p}
\end{corollary}

\begin{proof}
The subgraph of the graph $H_4$ induced by vertices in $A \cup B$ is isomorphic to $H_1$, which is not a PCG. Hence, by Lemma \ref{fl4}, the graph $H_4$ is also not a PCG. 
\qed \end{proof}

We now show that the graph $H_4$ belongs to the class $\mathcal{G}^4$ by showing that its complement graph $H_4^c$ has only one chordless cycle.

\begin{lemma}
    The graph $H_4$ belongs to the class $\mathcal{G}^4$.  
    
    \label{fl4c}
\end{lemma}

\begin{proof}
    In the graph $H_4$, from every vertex $v$ in the set $A \cup B$, there exists an edge to every vertex in $C$. Hence, in the graph $H_4^c$, from no vertex $v$ in the set $A \cup B$, there exists an edge to any vertex in $C$. Hence, in the graph $H_4^c$, the vertices in $C$ and the vertices in $A \cup B$  be in different connected components. 
    
    The connected component of the graph $H_4^c$ consisting of vertices in the set C is the complement of $C_5^c$, i. e., $C_5$, a chordless cycle. 
    
    The subgraph of the graph $H_4$ induced by the vertices in the set $A \cup B$ is isomorphic to the graph $H_1$. Therefore, the connected component(s) of the graph $H_4^c$ consisting of vertices in $A \cup B$ become the complement of $H_1$, i. e., $H_1^c$. Since $H_1 \in \mathcal{G}^1$, there is no chordless cycle in $H_1^c$.
    
    Hence, in total, the graph $H_4^c$ has only one chordless cycle and the proof is completed. 
\qed \end{proof}
From Corollary \ref{fl4p} and Lemma \ref{fl4c}, we get the following theorem. 

\begin{theorem}
    The class $\mathcal{G}^4$ contains some graphs which are not PCGs.
    \endproof
\end{theorem}

\section{Conclusion}
\label{sec:conclusion}
In this paper, we have proven that grid graphs are PCGs by giving a linear time algorithm to construct the PCT of a given grid graph. Since we have given a caterpillar construction for the PCT of a grid graph, this construction may help to find a complete characterization of PCGs admitting caterpillars as their witness trees. However, there are different superclasses of grid graphs like partial grid graphs, 3D grid graphs, king graphs etc and we do not know whether they are PCGs or not. An interesting future work might be to find all the properties of the PCTs of triangle free PCGs. Moreover, we have proven that two of the graph classes intermediate to the classes defined by the necessary condition and the sufficient condition are not PCGs leading to a conclusion that neither of the four intermediate graph classes defined by Hossain et al. \cite{JGAA419} are PCGs. These proofs may be helpful for a complete characterization of PCGs and finding the complexity class of the PCG recognition problem. 

\section*{Acknowledgement}

The authors would like to thank anonymous reviewers for
their suggestions for improving the presentation of the
paper. This work is supported by Basic Research Grant of BUET. 

%
%
%
 \bibliographystyle{elsarticle-num} 
 \bibliography{doc1}
%





\end{document}